
\magnification1200

\def\DJ{\leavevmode\setbox0=\hbox{D}\kern0pt\rlap
 {\kern.04em\raise.188\ht0\hbox{-}}D}

\def\hf{{\textstyle{1\over2}}}

\def\d{{\,\rm d}}
\def\e{\varepsilon}
\def\f{\varphi}
\def\G{\Gamma}

\def\s{\sigma}

\def\={\;=\;}
\def\zx{\zeta(\hf+ix)}
\def\zt{\zeta(\hf+it)}

\def\no{\noindent}  \def\kk{\kappa_h}
\def\R{\Re{\rm e}\,} \def\I{\Im{\rm m}\,} \def\s{\sigma}
\def\z{\zeta}

\def\no{\noindent} 

  \def\={\,=\,}
\def\hf{{\textstyle{1\over2}}}

\def\f{\varphi}
\def\Z{{\cal Z}}
\footnote{{}}{\hskip -2 mm 2001 {\it Mathematics Subject Classification}.
11M06}
\font\teneufm=eufm10
\font\seveneufm=eufm7
\font\fiveeufm=eufm5
\newfam\eufmfam
\textfont\eufmfam=\teneufm
\scriptfont\eufmfam=\seveneufm
\scriptscriptfont\eufmfam=\fiveeufm
\def\mathfrak#1{{\fam\eufmfam\relax#1}}

\font\tenmsb=msbm10
\font\sevenmsb=msbm7
\font\fivemsb=msbm5
\newfam\msbfam
      \textfont\msbfam=\tenmsb
      \scriptfont\msbfam=\sevenmsb
      \scriptscriptfont\msbfam=\fivemsb
\def\Bbb#1{{\fam\msbfam #1}}

\def \NN {\Bbb N}
\def \CC {\Bbb C}

  \def\rightheadline{{\hfil{\sevenrm The Mellin transform of the square of Riemann's
zeta-function
   }\hfil\tenrm\folio}}

  \def\leftheadline{{\tenrm\folio\hfil{\sevenrm
  Aleksandar Ivi\'c }\hfil}}
  \def\emptyheadline{\hfil}
  \headline{\ifnum\pageno=1 \emptyheadline\else
  \ifodd\pageno \rightheadline \else \leftheadline\fi\fi}

\font\bb=cmcsc10
\font\cc=cmcsc10 at 8pt

\font\kk=cmcsc10 at 12pt

\centerline{\kk The Mellin transform of the square of Riemann's
zeta-function   }

\medskip
\centerline{\bb  Aleksandar Ivi\'c}

\bigskip
{\bf Abstract.} Let $\Z_1(s) = \int_1^\infty |\zx|^2x^{-s}\d x\;
(\s = \R s > 1)$. A result concerning analytic continuation of
$Z_1(s)$ to $\CC$ is proved, and also a result relating the
order of $\Z_1(\s+it)\;(\hf \le \s \le 1, \,t\ge t_0)$
to the order of $\Z_1(\hf+it)$.

\bigskip

\centerline{\bf 1. Introduction}
\bigskip

Let $\Z_k(s)$, the (modified) Mellin transform of $|\zx|^{2k}$,
denote the analytic continuation of the function defined initially by
$$
\Z_k(s) \;:=\; \int_1^\infty |\zx|^{2k}x^{-s}\d x
\qquad(k\in\NN,\;\s = \R s > c(k) \;(>1)).
\leqno(1.1)
$$
This function, when $k=2$, was introduced by Y. Motohashi [15]
(see also [16]),
and its properties were further studied in [10] and [11].
The latter work also contains
some results on the function $\Z_1(s)$, which is the principal object of the
study in this paper. It was shown that $\Z_1(s)$ is regular for $\s > -3/4$,
except for a double pole at $s=1$. The principal part of the Laurent
expansion of $\Z_1(s)$ at $s=1$ is
$$
{1\over(s-1)^2} + {2\gamma - \log(2\pi)\over s - 1},\leqno(1.2)
$$
where $\gamma = -\G'(1) = 0.577215\ldots\,$ is Euler's constant.
M. Jutila [13] continued the study of $\Z_1(s)$ and proved
that $\Z_1(s)$ continues meromorphically to $\CC$, having only a double pole
at $s=1$ and at most double poles for $s = -1,-2,\ldots\,$.

\medskip
The object of this note is to prove some new results on $\Z_1(s)$. First we
make more precise Jutila's result on the analytic continuation of $\Z_1(s)$.
We have the following

\bigskip
{THEOREM 1}. {\it The function $\Z_1(s)$
continues meromorphically to $\CC$, having only a double pole
at $s=1$, and at most simple poles at $s = -1,-3,\ldots\,$. The principal
part of its Laurent expansion at $s=1$ is given by} (1.2).

\bigskip

Our second aim is to prove an order result for $\Z_1(s)$. This is
a result M. Jutila  mentioned in [13], but the term that I initially claimed,
namely $t^{1-2\s+\e}$ is too optimistic. It appears that what can
be proved is contained in

\bigskip
THEOREM 2. {\it We have, for $\hf \le
\s \le 1,\;t\ge t_0 > 0$,}
$$
\Z_1(\s + it) \;\ll_\e\; t^{{1\over2}-\s+\e}
\max_{t-t^\e\le v\le t+t^\e}|\Z_1(\hf+iv)| + (t^{9-16\s\over7}+t^{-1})\log t.
\leqno(1.3)
$$

\bigskip
{\bf Corollary} (see M. Jutila [13]). For $\hf \le \s \le 1,\;t\ge t_0$
we have
$$
\Z_1(\s+it) \ll_\e t^{{5\over6}-\s+\e}.\leqno(1.4)
$$

\medskip
The bound $\Z_1(\hf+it) \ll_\e t^{1/3+\e}$ was mentioned in [11], and
its proof was elaborated by M. Jutila in [13]. If one inserts this
bound in (1.3), then (1.4) follows immediately.
We note that here and later $\e$ denotes arbitrarily small constants,
not necessarily the same ones at each occurrence.
It  can be also proved that, for any given $\e > 0$,
$$
\int_1^T|\Z_k(\s + it)|^2\d t \;\gg_\e\; T^{2-2\s-\e}\qquad(k = 1,2;\;
\hf < \s < 1).\leqno(1.5)
$$
The bound (1.5) for $k = 2$  appeared in my paper [10], and the
proof of the bound when $k=1$ is on similar lines, so that the details
will not be given here. It is plausible that
$$
\max_{t-t^\e\le v\le t+t^\e}|\Z_1(\hf+iv)| \gg 1\leqno(1.6)
$$
holds, but  at present I am unable to prove (1.6).

\bigskip
\centerline{2. \bf Proof of Theorem 1}

\bigskip
Let
$$
L_k(s) := \int_0^\infty |\zx|^{2k}{\rm e}^{-sx}\d x\qquad
(\s = \R s > 0,\,k\in\NN)\leqno(2.1)
$$
be the Laplace transform of $|\zx|^{2k}$. We are interested in
using the expression for $L_1(s)$ when $s = 1/T,\,T\to\infty$.
Such a formula has been known for a long time and is due to H.
Kober [14] (see also a proof in [17, Chapter 9]). The functions
$L_k(s)$ when $k = 1,2$ were studied by F.V. Atkinson [1], [2].
The function $L_2(s)$ was considered by the author [8], [9], and
the approach to mean value of Dirichlet series by Laplace
transforms by M. Jutila [12]. Kober's result on $L_1(s)$ is that,
for any integer $N\ge0$,
$$
L_1(\s) = {\gamma-\log(2\pi\s)\over2\sin{\sigma\over2}}
+ \sum_{n=0}^Nd_n\s^n + O_N(\s^{N+1})\quad(\s\to0+),\leqno(2.2)
$$
where $d_n$ are suitable constants. Note that for
$\s = 1/T,\,T\to\infty$ we have
$$
\eqalign{&
{\gamma-\log\left({2\pi\over T}\right)\over2\sin\left({1\over2T}\right)}
= {\log\left({T\over2\pi}\right) + \gamma\over{1\over T} -
{1\over24T^3} + {1\over16\cdot5!T^5} - \cdots}\cr&
= \left(\log\left({T\over2\pi}\right) +
\gamma\right)\sum_{n=0}^\infty c_nT^{1-2n}\cr}
$$
with the coefficients $c_n$ that may be explicitly evaluated. Therefore
for any integer $N\ge0$
$$
L_1\left({1\over T}\right) =
\left(\log\left({T\over2\pi}\right) +
\gamma\right)\sum_{n=0}^N a_nT^{1-2n}
+ \sum_{n=0}^N b_nT^{-2n} + O_N(T^{-1-2N}\log T)\leqno(2.3)
$$
with the coefficients $a_n, b_n$ that may be explicitly evaluated.
In particular, we have
$$
a_0 \=\ 1,\quad b_0\=\ \pi.\leqno(2.4)
$$
It is clear from Kober's formula (2.2) that $a_0 = 1$.
To see that $b_0 = \pi$
one can use the work of Hafner--Ivi\'c [4] (this was mentioned
by Conrey et al. [3]),
which will be stated now in detail, as it will be needed also for the proof
of Theorem 2. Let, as usual, for $T\ge0$,
$$
E(T) \;=\; \int_0^T|\zt|^2\d t - T\left(\log\left({T\over2\pi}\right) +
2\gamma -1\right) \quad (E(0) = 0)
$$
denote the error term in the mean square formula for $|\zt|$.
Let further
$$
G(T) := \int_0^T E(t)\d t - \pi T,\quad G_1(T) := \int_0^T G(t)\d t.
\leqno(2.5)
$$
Then  Hafner--Ivi\'c [4] proved (see also [6])
$$
G(T) \= S_1(T;N) - S_2(T;N) + O(T^{1/4})\leqno(2.6)
$$
with
$$
S_1(T;N) = 2^{-3/2}\sum_{n\le N}(-1)^nd(n)n^{-1/2}\left({\rm arsinh}\,
\sqrt{{\pi n
\over2T}}\,\right)^{-2}\left({T\over2\pi n}+{1\over4}\right)^{-1/4}
\sin(f(T,n)),
$$
$$
S_2(T;N) = \sum_{n\le N'}d(n)n^{-1/2}\left(\log {T\over2\pi n}\right)^{-2}
\sin(g(T,n)),
$$
$$
f(T,n) \= 2T\,{\rm arsinh}\,\sqrt{{\pi n\over2T}} +
\sqrt{2\pi nT + \pi^2n^2} - {\pi\over4},
$$
$$
g(T,n) \= T\log\left({T\over2\pi n}\right) - T + {\pi\over4},
\;{\rm arsinh}\,x = \log(x + \sqrt{x^2+1}),
$$
$$
AT < N < A'T\;(0 < A < A' \,{\rm constants}),\;
N' = {T\over2\pi} + {N\over2} - \sqrt{{N^2\over4} + {NT\over2\pi}}.
$$
We use Taylor's formula (see [7, Lemma 3.2])
to simplify  (2.6). Then we obtain
$$
G(T) = 2^{-1/4}\pi^{-3/4}T^{3/4}\sum_{n=1}^\infty
(-1)^nd(n)n^{-5/4}\sin(\sqrt{8\pi nT} - {\pi\over4}) + O(T^{2/3}\log T),
$$
so that
$$
G(T) \= O(T^{3/4}),\qquad G(T) \= \Omega_\pm(T^{3/4}).\leqno(2.7)
$$
We use (2.6) with $T = t, N = T, T \le t \le 2T$ and apply the first
derivative test ([5, Lemma 2.1]) to deduce that
$$
\int_T^{2T}G(t)\d t \;\ll\; T^{5/4},
$$
since the term $O(1)$ in [5, eq. (15.29)] in Atkinson's formula is in fact
$O(1/T)$. Hence
$$
G_1(T) \;=\; \int_0^T G(t)\d t \;\ll\;T^{5/4}.\leqno(2.8)
$$

\medskip
To prove that $b_0=\pi$ (cf. (2.4)) we note that,
using the definition of $E(T)$ and applying integration by parts,
$$
\eqalign{&
L_1\left({1\over T}\right) = \int_0^\infty {\rm e}^{-t/T}
\left(\log\left({t\over2\pi}\right) + 2\gamma  + E'(t)\right)\d t\cr&
= T\int_0^\infty {\rm e}^{-x}(\log x + \log {T\over2\pi} + 2\gamma)\d x
+ {1\over T}\int_0^\infty E(t){\rm e}^{-t/T}\d t\cr&
= T\left(\log\left({T\over2\pi}\right) + \gamma\right) +
{1\over T^2}\int_0^\infty \int_0^tE(u)\d u\cdot {\rm e}^{-t/T}\d t\cr&
= T\left(\log\left({T\over2\pi}\right) + \gamma\right) +
{1\over T^2}\int_0^\infty (\pi t + O(t^{3/4})){\rm e}^{-t/T}\d t\cr&
= T\left(\log\left({T\over2\pi}\right) + \gamma\right) + \pi + O(T^{-1/4}),
\cr}\leqno(2.9)
$$
where we also used the $O$-bound of (2.7) and
$$
\int_0^\infty {\rm e}^{-x}\log x\d x = \G'(1) = -\gamma.
$$
A comparison of (2.3)  and (2.9) proves that $b_0 = \pi$, as asserted
by (2.4).

\medskip
We return now to the proof of Theorem 1.
Let
$$
{\bar L}_k(s) := \int_1^\infty|\z(\hf+iy)|^{2k}{\rm e}^{-ys}
\d y\quad(k\in\NN,\; \s = \R s > 0). \leqno(2.10)
$$
Then we have by absolute convergence, taking $\s $ sufficiently large
and making the change of variable $xy = t$,
$$\eqalign{&
\int_0^\infty {\bar L}_k(x) x^{s-1}\d x = \int_0^\infty
\left(\int_1^\infty|\z(\hf + iy)|^{2k}{\rm e}^{-yx}\d y\right)
x^{s-1}\d x\cr&
= \int_1^\infty|\z(\hf + iy)|^{2k}\left(\int_0^\infty
x^{s-1}{\rm e}^{-xy}\d x\right)\d y\cr&
= \int_1^\infty|\z(\hf + iy)|^{2k}y^{-s}\d y \int_0^\infty
{\rm e}^{-t}t^{s-1}\d t = {\cal Z}_k(s)\G(s).\cr}\leqno(2.11)
$$
Further we have
$$\eqalign{&
\int_0^\infty {\bar L}_1(x) x^{s-1}\d x = \int_0^1{\bar L}_1(x) x^{s-1}\d x
+ \int_1^\infty{\bar L}_1(x) x^{s-1}\d x\cr&
= \int_1^\infty{\bar L}_1(1/x) x^{-1-s}\d x + A(s)\quad(\s>1),\cr}
$$
where
$$
A(s) \;:=\; \int_1^\infty{\bar L}_1(x) x^{s-1}\d x
$$
is an entire function. Since
$$
{\bar L}_1(1/x) =  L_1(1/x) - \int_0^1|\z(\hf+iy)|^2{\rm e}^{-y/x}
\d y\qquad(x\ge1),
$$
it follows from (2.11) with $k=1$ by analytic continuation that, for $\s>1$,
$$\eqalign{
{\cal Z}_1(s)\G(s) &= \int_1^\infty L_1(1/x) x^{-1-s}\d x
- \int_1^\infty\bigl(\int_0^1|\z(\hf+iy)|^2{\rm e}^{-y/x}
\d y\bigr)x^{-1-s}\d x + A(s)\cr&
= I_1(s) - I_2(s) + A(s),\cr}\leqno(2.12)
$$
say. Clearly for any integer $M\ge1$
$$\eqalign{
I_2(s)&= \int_1^\infty\int_0^1|\z(\hf+iy)|^2
\left(\sum_{m=0}^M{(-1)^m\over m!}\left({y\over x}\right)^m +
O_M(x^{-M-1})\right)
\d y \,x^{-1-s}\d x\cr&
= \sum_{m=0}^M{(-1)^m\over m!}h_m\cdot{1\over m+s} + H_M(s),\cr}\leqno(2.13)
$$
say, where $H_M(s)$ is a regular function of $s$ for $\s > -M-1$, and
$$
h_m \;:=\;\int_0^1|\z(\hf+iy)|^2y^m\d y
$$
is a constant. Inserting (2.3) in $I_1(s)$ in (2.12) we have, for $\s>1$,
$$\eqalign{
I_1(s) &= \int_1^\infty (\log {x\over2\pi}+\gamma)
\sum_{n=0}^Na_nx^{-2n-s}\d x + \int_1^\infty \sum_{n=0}^Nb_nx^{-1-2n-s}\d x
+ K_N(s)\cr&
= \sum_{n=0}^Na_n\left({1\over(2n+s-1)^2} +
{\gamma-\log2\pi\over 2n+s-1}\right)
+  K_N(s),\cr}\leqno(2.14)
$$
say, where $K_N(s)$ is regular for $\s > -2N$. Taking $M = 2N$
it follows from (2.12)--(2.14) that
$$\eqalign{
\Z_1(s)\G(s) &= \sum_{n=0}^Na_n\left({1\over(2n+s-1)^2} +
{\gamma-\log2\pi\over 2n+s-1}\right)\cr&
+ \sum_{m=0}^{2N}{(-1)^m\over m!}h_m\cdot{1\over m+s} + R_N(s),\cr}
\leqno(2.15)
$$
say, where $R_N(s)$ is a regular function of $s$ for $\s > -2N$. This holds
initially for $\s >1$, but by analytic continuation it holds for $\s > -2N$.
Since $N$ is arbitrary
and $\G(s)$ has no zeros, it follows that (2.15)
provides meromorphic continuation
of $\Z_1(s)$ to $\CC$. Taking into
account that $\G(s)$ has simple poles at $s = -m\;(m=0,1,2,\ldots\,)$ with
residues $(-1)^m\over m!$, and that near $s=1$ we have the Taylor expansion
$$
{1\over \G(s)} = 1 + \gamma(s-1) + \sum_{n=2}^\infty {f_n\over n!}(s-1)^n
$$
with $f_n = (1/\G(s))^{(n)}\Big|_{s=1}$, it follows that
$$\eqalign{
\Z_1(s) &=
{1\over(s-1)^2} + {2\gamma - \log(2\pi)\over s - 1}\cr&
+ {1\over\G(s)}\left\{\sum_{n=1}^Na_n\left({1\over(2n+s-1)^2} +
{\gamma-\log2\pi\over 2n+s-1}\right)\right\}\cr&
+ {1\over\G(s)}\left(\sum_{m=0}^{2N}{(-1)^m\over m!}h_m\cdot{1\over m+s}
\right) + U_N(s),\cr}      \leqno(2.16)
$$
say, where $U_N(s)$ is a regular function of $s$ for $\s > -2N$.
Moreover, the function  $1/(\G(s)(2n+s-1))$  is regular
for $n \in\NN$ and $s\in\CC$. Therefore (2.16) provides analytic
continuation to $\CC$, showing that besides $s=1$
the only poles of $\Z_1(s)$ can be
simple poles at $s = 1-2n$ for $n \in\NN$, as asserted by Theorem 1.
With more care the residues at these poles could be explicitly evaluated.
\bigskip

\centerline
{\bf 3. Proof Theorem 2}

\bigskip

 We suppose that $X\gg t$ and start from
$$
{\cal Z}_1(s) = \left(\int_1^X + \int_X^\infty\right)|\zx|^2x^{-s}\d x
\qquad(\s = \R s > 1).
$$
We use $|\zx|^2 = \log(x/2\pi) + 2\gamma + (E(x)-\pi)'$ and
integrate by parts to obtain
$$\eqalign{
{\cal Z}_1(s) &= \int_1^{t^{1-\e}}|\zx|^2x^{-s}\d x +
\int_{t^{1-\e}}^X|\zx|^2x^{-s}\d x \cr& +
{X^{1-s}\over s - 1}\left({1\over s-1} + \log X + 2\gamma - \log2\pi\right)
 \cr& -(E(X)-\pi)X^{-s}   + s\int_X^\infty (E(x)-\pi)x^{-s-1}\d x \cr&
= I_1 + I_2 + O(X^{-\s}+t^{-1}X^{1-\s}\log X) + I_3,
\cr}\leqno(3.1)
$$
provided $E(X) = 0$. But as proved in [6], every interval
$[T, T + C\sqrt{T}]$ with sufficiently large $C>0$ contains
a zero of $E(t)$, hence we can assume that $E(X) = 0$
is fulfilled with suitable $X$.
Since $E(x)$ is $\asymp x^{1/4}$ in mean square (see e.g., [5] or [7]),
it follows by analytic continuation
that (3.1) is valid for $\s > {1\over4},\,t \ge t_0> 0,\,
 t \le X \le t^A\,(A \ge 1).$
By repeated integration by parts it is found that $I_1 \ll t^{-1}$.
The integral $I_2$ is split in $O(\log T)$ subintegrals of the form
$$
J_K \= \int_{K/2}^{5K/2}\f(x)|\zx|^2x^{-s}\d x,
$$
where $\f(x) \in C^\infty$ is a nonnegative, smooth function supported
in $[K/2,\,5K/2]\,$ that is equal
to unity in $[K,\,2K]\,$. To bound $J_K$ one can start
from (1.1) with $k=1$ and use the Mellin inversion
formula (see [10] for a detailed discussion concerning the
analogous situation with the fourth moment of $|\zt|$), which yields
$$
|\zx|^2 \= {1\over2\pi i}\int_{(1+\e)}{\cal Z}_1(s)x^{s-1}\d s
\qquad(x > 1).
$$
Here we replace the line of integration  by the contour ${\cal L}$,
consisting of the same straight line from which the segment
$[2+\e-i,\,1+\e+i]$ is removed and replaced by a circular arc
of unit radius, lying
to the left of the line, which passes over the pole $s =1 $ of
the integrand. By the residue theorem we have
$$
|\zx|^2 \= {1\over2\pi i}\int_{\cal L}{\cal Z}_1(s)x^{s-1}\d s
+ \log \left({x\over2\pi}\right) + 2\gamma  \qquad(x > 1).
$$
This gives
$$
J_K \= {1\over2\pi i}\int_{{\cal L}}{\cal Z}_1(w)
\int_{K/2}^{5K/2}\f(x)x^{w-s-1}\d x\d w  + O(t^{-1}).
$$
Repeated integration by parts in the $x$--integral shows that only
the portion $|\I w - \I s| \le t^\e$ gives a non-negligible
contribution. Hence replacing ${\cal L}$ by the line $\R w = \hf$
we obtain
$$
J_K \;\ll_\e\; K^{{1\over2}-\s}t^\e\max_{t-t^\e\le v\le t+t^\e}|\Z_1(\hf+iv)|
+ t^{-1}.\leqno(3.2)
$$
To  estimate $I_3$ in (3.1) we integrate by parts and use (2.7) and (2.8).
We obtain
$$
\eqalign{&
\int_X^\infty (E(x)- \pi)x^{-s-1}\d x = -G(X)X^{-s-1}+(s+1)\int_X^\infty
G(x)x^{-s-2}\d x\cr&
= O(X^{-1/4-\s}) - (s+1)G_1(X)X^{-s-2} + (s+1)(s+2)\int_X^\infty
G_1(x)x^{-s-3}\d x\cr&
\ll X^{-1/4-\s} + t^2X^{-3/4-\s}.\cr}
$$
Thus it follows from (3.1) and (3.2) that we have, for
$\hf \le \s \le 1,\;t\ge t_0$,
$$\eqalign{
\Z_1(\s + it) &\ll_\e t^{{1\over2}-\s+\e}
\max_{t-t^\e\le v\le t+t^\e}|\Z_1(\hf+iv)| +
(X^{1-\s}t^{-1}+ t^{-1})\log t  \cr&
+ X^{-\s} + tX^{-1/4-\s} + t^3X^{-3/4-\s}.\cr}\leqno(3.3)
$$
In (3.3) we choose $X$ to satisfy
$$
X^{1-\s}t^{-1} = t^3X^{-3/4-\s}, \quad{\rm i.e.,}\quad X = t^{16/7}.
$$
Then $tX^{-1/4-\s}+ X^{-\s} \ll t^3X^{-3/4-\s}$, and (1.3) follows from (3.3).

\bigskip
In conclusion, one may ask what should be the true order of
magnitude of $\Z_1(\s + it) $. This is a difficult question, but
it seems reasonable to expect that
$$
\Z_1(\s + it) \ll_\e t^{{1\over2}-\s+2\mu({1\over2})+\e}\qquad(\hf \le
\s \le 1,\;t\ge t_0 > 0)\leqno(3.4)
$$
holds, where
$$
\mu(\s) := \limsup_{t\to \infty}\,{\log|\z(\s+it)|\over\log t}
$$
is the Lindel\"of exponent of $\zeta(s)$. One could try to obtain
(3.4) from (1.3) by refining the estimation of ${\cal Z}_1(s)$ in
[11]. This procedure leads to an exponential sum with the divisor
function which is estimated as $\ll_\e t^{2\mu({1\over2})+\e}$
(this, of course, is non-trivial and needs elaboration).

\bigskip
\centerline{\kk References}
\bigskip

\item{[1]} F.V. Atkinson, The mean value of the zeta-function on
the critical line, {\it Quart. J. Math. Oxford} {\bf 10}(1939), 122-128.

\item {[2]} F.V. Atkinson, The mean value of the zeta-function on
the critical line, {\it Proc. London Math. Soc.} {\bf 47}(1941), 174-200.

\item{[3]} J.B. Conrey, D.W. Farmer, J.P. Keating, M.O. Rubinstein and N.C. Snaith,
Integral moments of $L$-functions, 71pp, {\tt ArXiv.math.NT/0206018}.

\item {[4]}  J.L. Hafner and A. Ivi\'c,
On the mean square of the Riemann zeta-function
on the critical line, {\it J. Number Theory } {\bf 32}(1989), 151-191.

\item {[5]} A. Ivi\'c,  The Riemann zeta-function, {\it John Wiley
and Sons}, New York, 1985.

\item {[6]} A. Ivi\'c, Large values of certain number-theoretic error terms,
{\it Acta  Arith.} {\bf56}(1990), 135-159.

\item {[7]} A. Ivi\'c,  Mean values of the Riemann zeta-function,
LN's {\bf 82}, {\it Tata Institute of Fundamental Research},
Bombay,  1991 (distr. by Springer Verlag, Berlin etc.).

\item{[8]}  A. Ivi\'c,  On the fourth moment of the Riemann
zeta-function, {\it Publs. Inst. Math. (Belgrade)} {\bf 57(71)}
(1995), 101-110.

\item{[9]} A. Ivi\'c, The Laplace transform of the fourth
moment of the zeta-function,
{\it Univ. Beograd. Publ. Elektrotehn. Fak. Ser. Mat.} {\bf11}(2000), 41-48.

\item{[10]} A. Ivi\'c, On some conjectures and results for the
Riemann zeta-function, {\it Acta. Arith.} {\bf99}(2001), 115-145.

\item {[11]} A. Ivi\'c, M. Jutila and Y. Motohashi, The Mellin
transform of powers of  the Riemann zeta-function, {\it Acta Arith.}
{\bf95}(2000), 305-342.

\item{[12]}  M. Jutila, Mean values of Dirichlet series via Laplace
transforms,
in {\it ``Analytic Number Theory"} (ed. Y. Motohashi), London Math. Soc.
 LNS {\bf247}, {\it Cambridge University Press}, Cambridge, 1997, 169-207.

\item {[13]}  M. Jutila, The Mellin transform of the square of Riemann's
zeta-function, {\it Periodica Math. Hung.} {\bf42}(2001), 179-190.

\item {[14]} H. Kober, Eine Mittelwertformel der Riemannschen Zetafunktion,
{\it Compositio Math.} {\bf3}(1936), 174-189.

\item {[15]} Y. Motohashi,  A relation  between the Riemann
zeta-function and the hyperbolic Laplacian, {\it Annali Scuola Norm.
Sup. Pisa, Cl. Sci. IV ser.} {\bf 22}(1995), 299-313.

\item{[16]} Y. Motohashi, Spectral  theory of the Riemann zeta-function,
Cambridge University Press, 1997.

\item {[17]} E.C. Titchmarsh,  The Theory of the Riemann
Zeta-Function (2nd ed.), {\it Clarendon  Press}, Oxford, 1986.


\vfill

  \no
\cc Aleksandar Ivi\'c

\no Katedra Matematike RGF-a

\no Universitet u Beogradu

\no \DJ u\v sina 7, 11000 Beograd

\no Serbia and Montenegro

\no e-mail: \sevenbf aivic@matf.bg.ac.yu,

\no {\hskip 10mm} ivic@rgf.bg.ac.yu

\bye